\documentclass[12pt,thmsa]{article}%
\usepackage{amssymb}
\usepackage{sw20bams}
\usepackage{amsmath}
\usepackage{amsfonts}
\usepackage{graphicx}%
\providecommand{\U}[1]{\protect\rule{.1in}{.1in}}
\begin{document}

\author{Steven Finch}
\title{In Limbo:\ Three Triangle Centers }
\date{June 3, 2014}
\maketitle

\begin{abstract}
Yet more candidates are proposed for inclusion in the Encyclopedia of Triangle
Centers. \ Our focus is entirely on simple calculations.

\end{abstract}

\footnotetext{Copyright \copyright \ 2014 by Steven R. Finch. All rights
reserved.}The best-known triangle centers:

\begin{itemize}
\item incenter (intersection of three angle bisectors)

\item centroid (intersection of three medians)

\item circumcenter (intersection of three perpendicular bisectors)

\item orthocenter (intersection of three altitudes)
\end{itemize}

\noindent are the first four listed in Kimberling's famous database \cite{Kb}.
\ Thousands more appear. \ A\ recent addition is the electrostatic center
\cite{AK}. \ Rigorous definition of a triangle center involves a real function
$f$ defined on the set of all possible triples $(a,b,c)$ of triangle sides and
satisfying certain properties. \ We forego such requirements (hence the phrase
\textquotedblleft in limbo\textquotedblright) and informally propose three
more triangle centers:

\begin{itemize}
\item equiareal disk center (associated with Fraenkel asymmetry \cite{FMP, Fi})

\item illuminating center (also called a Shibata streetlight \cite{Sh, OH})

\item thermodynamic center (also called a \textquotedblleft hot
spot\textquotedblright\ \cite{Pa, BMS})
\end{itemize}

\noindent in the hope that someone else can pick up where we leave off. \ The
latter notion, like the electrostatic center, has its origins in physics.

Consider three triangles $T_{1}$, $T_{2}$, $T_{3}$ with vertices%
\[%
\begin{array}
[c]{ccc}%
\left\{  0,0\right\}  ,\left\{  1,0\right\}  ,\left\{  0,1\right\}  , &  &
\text{isosceles right triangle;}%
\end{array}
\]%
\[%
\begin{array}
[c]{ccc}%
\left\{  0,0\right\}  ,\left\{  1,0\right\}  ,\left\{  0,\sqrt{3}\right\}  , &
& 30^{\circ}\text{-}60^{\circ}\text{-}90^{\circ}\text{ triangle;}%
\end{array}
\]%
\[%
\begin{array}
[c]{ccc}%
\left\{  0,0\right\}  ,\left\{  6,0\right\}  ,\left\{  -\frac{13}{3}%
,\frac{4\sqrt{35}}{3}\right\}  , &  & 6\text{-}9\text{-}13\text{ triangle.}%
\end{array}
\]
Our humble contribution is the calculation of triangle centers for these
cases. \ We make no claim of originality or special insight. \ If our paper
starts a conversation, leading perhaps to future inclusion of the three
centers in \cite{Kb}, then our efforts will be justified. \ 

\section{Equiareal Disk Center}

Given a triangle $T$, let $\left\vert T\right\vert $ denote its area. The
Fraenkel asymmetry of $T$ is defined to be\textbf{\ }%
\[
\alpha(T)=\inf\left\{  \frac{\left\vert (T\smallsetminus D)\cup
(D\smallsetminus T)\right\vert }{\left\vert T\right\vert }:D\text{ a disk with
}\left\vert D\right\vert =\left\vert T\right\vert \right\}  .
\]
The numerator contains the symmetric difference of $T$ and $D$. \ One could
imagine a similar definition involving disks having the same perimeter as $T$,
rather than area, but we leave this variation for other people to explore.
\ The infimum $\alpha(T)$ is achieved for some disk centered at a unique
interior point of $T$. Locating this point is a challenging exercise in
calculus. \ We illustrate the necessary partitioning of $T_{1}$, $T_{2}$,
$T_{3}$ in Figures 1, 2, 3 respectively. \ The color red is used for vertical
rectangles of width $dx$; green is used for horizontal rectangles of width
$dy$; black is used to further subdivide certain cells of the partition.
\ Details of other feasible configurations of the triangle and disk are
omitted for brevity's sake.%

%TCIMACRO{\FRAME{ftbpFU}{6.0243in}{5.9491in}{0pt}{\Qcb{Isosceles right triangle
%$T_{1}$}}{}{tri1.eps}{\special{ language "Scientific Word";  type "GRAPHIC";
%width 6.0243in;  height 5.9491in;  depth 0pt;  original-width 8.7398in;
%original-height 8.6265in;  cropleft "0";  croptop "1";  cropright "1";
%cropbottom "0";  filename 'tri1.eps';file-properties "XNPEU";}} }%
%BeginExpansion
\begin{figure}[ptb]%
\centering
\includegraphics[
height=5.9491in,
width=6.0243in
]%
{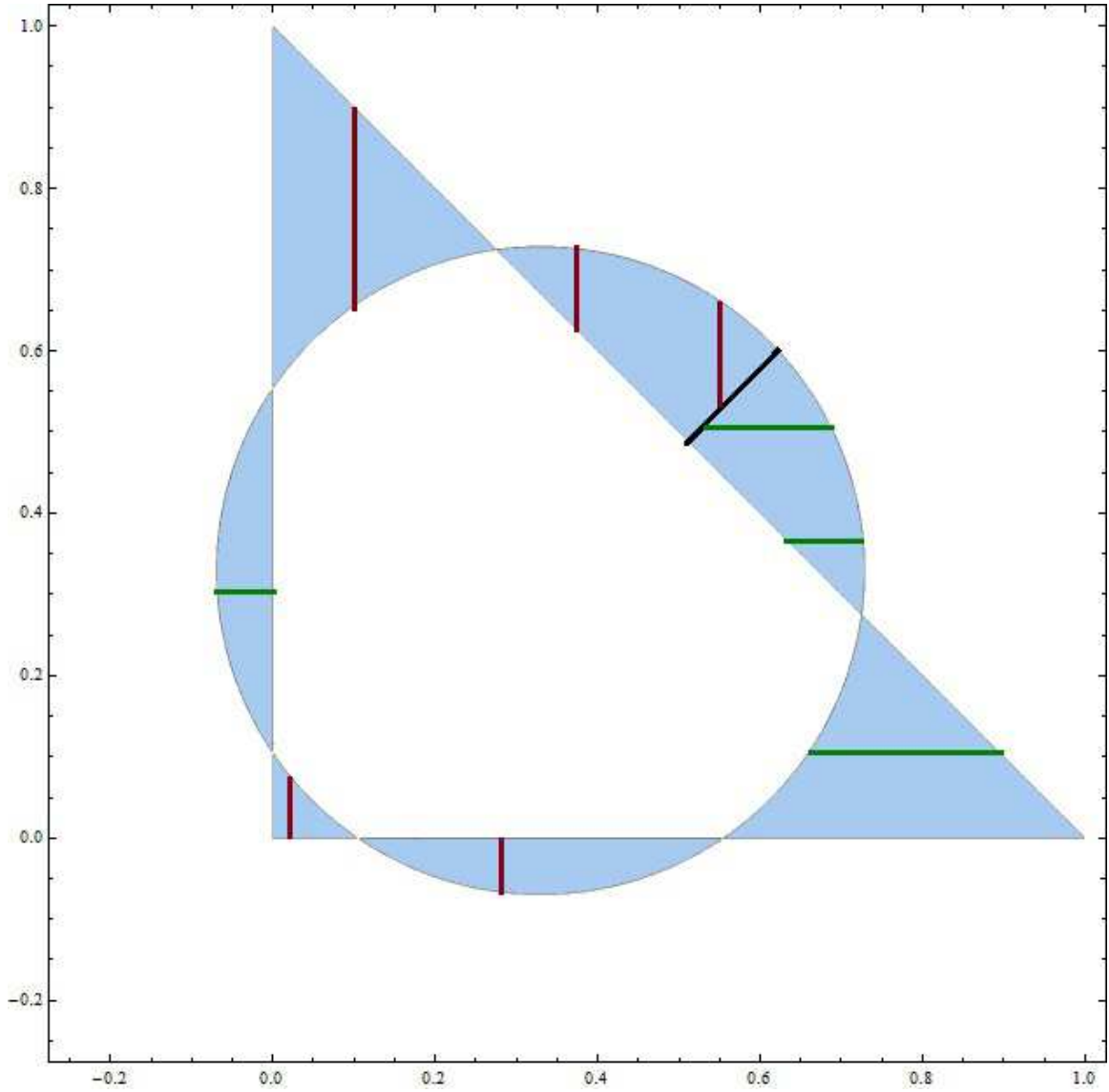}%
\caption{Isosceles right triangle $T_{1}$}%
\end{figure}
%EndExpansion

%

%TCIMACRO{\FRAME{ftbpFU}{6.1428in}{6.1134in}{0pt}{\Qcb{$30^{\circ}$-$60^{\circ
%}$-$90^{\circ}$ triangle $T_{2}$}}{}{tri2.eps}%
%width 6.1428in;  height 6.1134in;  depth 0pt;  original-width 3.627in;
%original-height 3.6097in;  cropleft "0";  croptop "1";  cropright "1";
%cropbottom "0";  filename 'tri2.eps';file-properties "XNPEU";}} }%
%BeginExpansion
\begin{figure}[ptb]%
\centering
\includegraphics[
height=6.1134in,
width=6.1428in
]%
{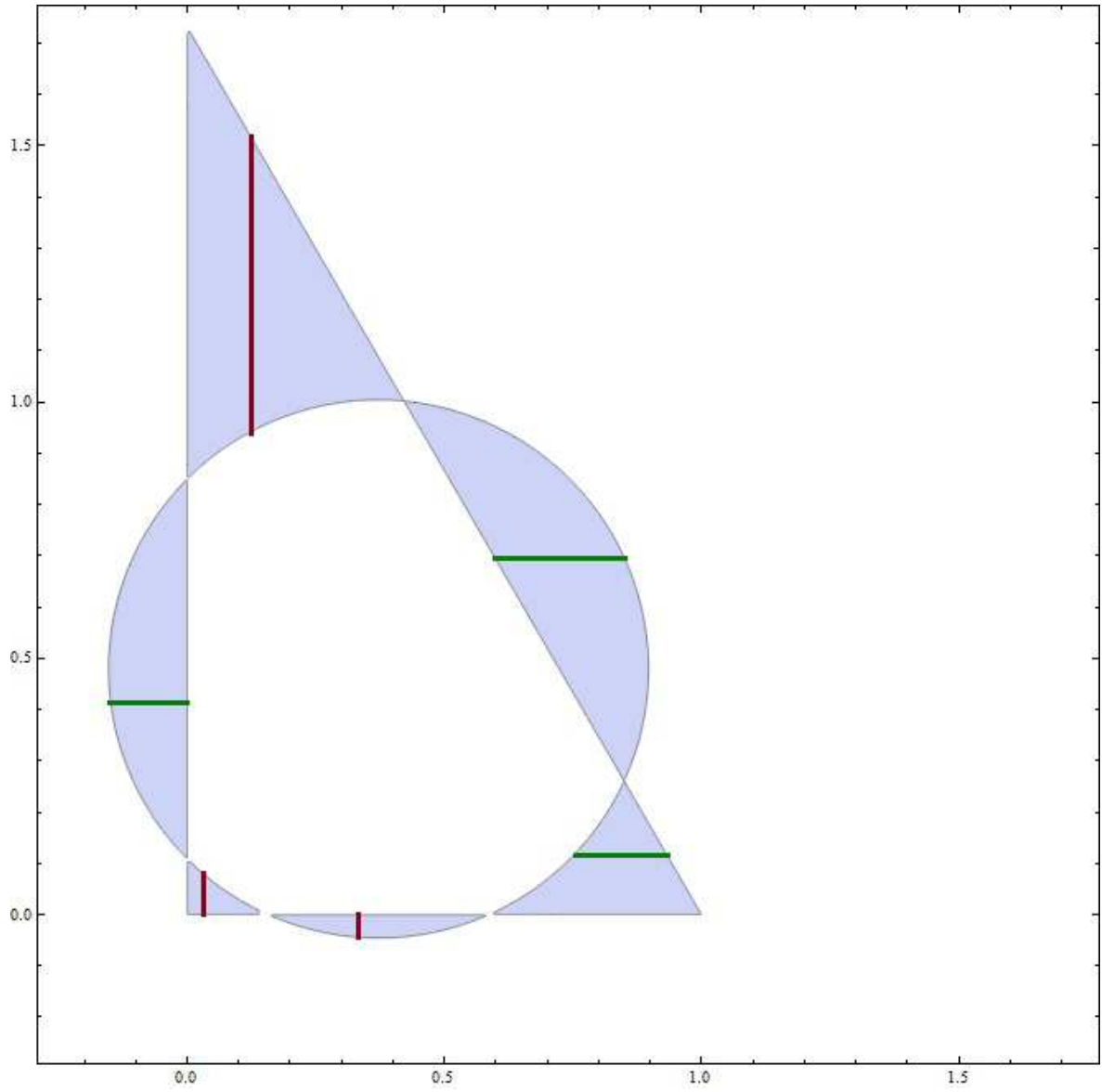}%
\caption{$30^{\circ}$-$60^{\circ}$-$90^{\circ}$ triangle $T_{2}$}%
\end{figure}
%EndExpansion
%TCIMACRO{\FRAME{ftbpFU}{5.9819in}{5.9819in}{0pt}{\Qcb{$6$-$9$-$13$ triangle
%$T_{3}$}}{}{tri3.eps}{\special{ language "Scientific Word";  type "GRAPHIC";
%width 5.9819in;  height 5.9819in;  depth 0pt;  original-width 3.2906in;
%original-height 3.2906in;  cropleft "0";  croptop "1";  cropright "1";
%cropbottom "0";  filename 'tri3.eps';file-properties "XNPEU";}} }%
%BeginExpansion
\begin{figure}[ptb]%
\centering
\includegraphics[
height=5.9819in,
width=5.9819in
]%
{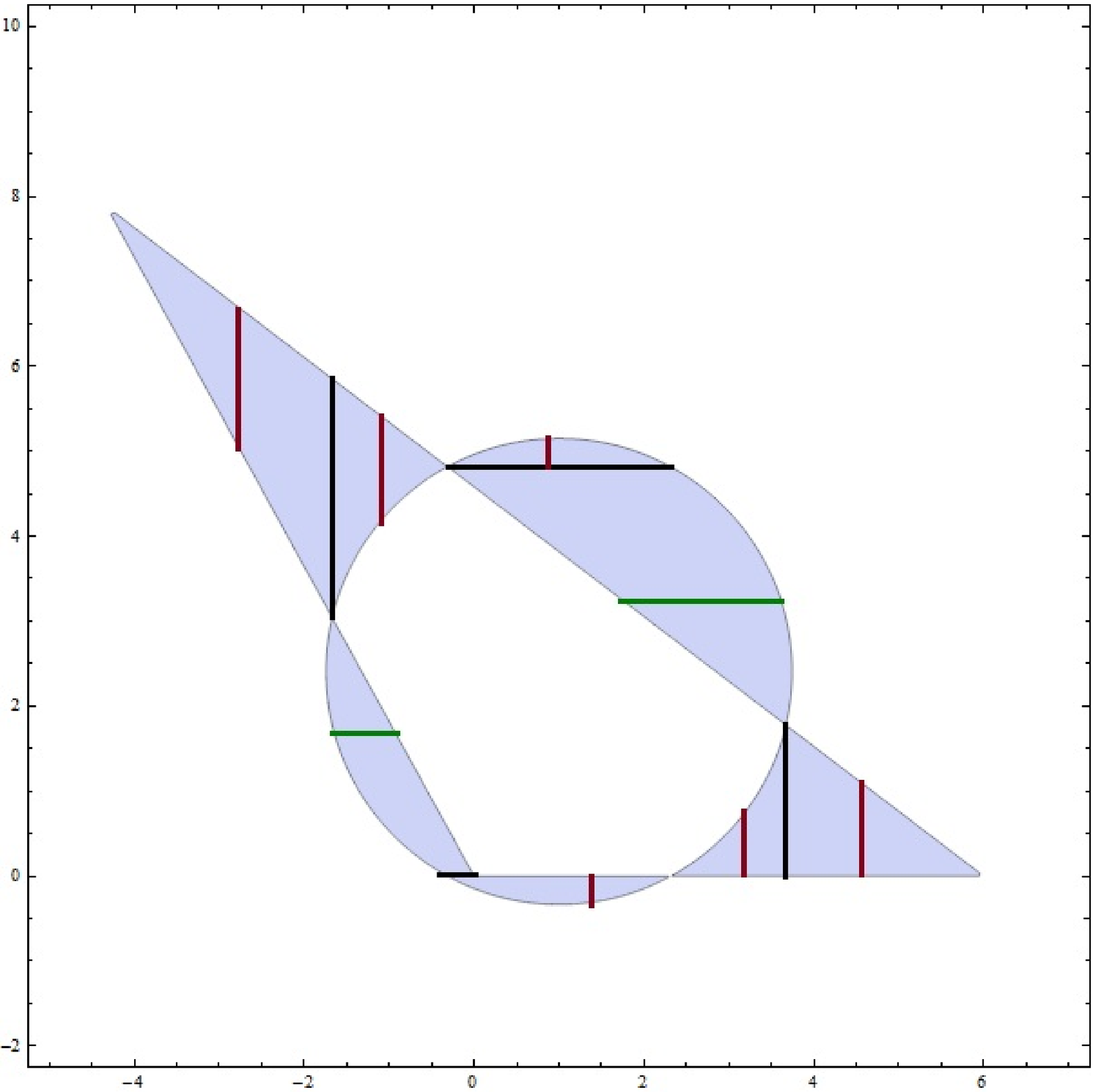}%
\caption{$6$-$9$-$13$ triangle $T_{3}$}%
\end{figure}
%EndExpansion

\subsection{Isosceles Right Triangle}

For $T_{1}$, the center must be on the diagonal line $y=x$ by symmetry.
\ Since $\left\vert T_{1}\right\vert =1/2$, the optimal circle $C$ has the
form%
\[
(x-t)^{2}+(y-t)^{2}=\dfrac{1}{2\pi}.
\]
The intersection of the line $x+y=1$ and $C$ yields two points $(p_{1},q_{1})
$, $(p_{2},q_{2})$ with $p_{1}<p_{2}$, $q_{1}>q_{2}$. \ The intersection of
the line $x=0$ and $C$ yields two points $(p_{3},q_{3})$, $(p_{4},q_{4})$ with
$p_{3}=p_{4}=0$, $q_{3}>q_{4}$. \ The intersection of the line $y=0$ and $C$
yields two points $(p_{5},q_{5})$, $(p_{6},q_{6})$ with $p_{5}<p_{6}$,
$q_{5}=q_{6}=0$. \ The intersection of the line $y=x$ and $C$ yields two
points; we select $(p_{7},q_{7})$ to be the point with the larger
$x$-coordinate. \ Everything can easily be made explicit, for example,%
\[%
\begin{array}
[c]{ccccc}%
p_{1}=\dfrac{1}{2}\left(  1-\sqrt{\dfrac{1}{\pi}-(1-2t)^{2}}\right)  , &  &
p_{5}=t-\sqrt{\dfrac{1}{2\pi}-t^{2}}, &  & p_{7}=\dfrac{1}{2\sqrt{\pi}}+t.
\end{array}
\]
The northwestern triangle has area\
\[%
%TCIMACRO{\dint \limits_{0}^{p_{1}}}%
%BeginExpansion
{\displaystyle\int\limits_{0}^{p_{1}}}
%EndExpansion
\left[  (1-x)-\left\{  t+\sqrt{\dfrac{1}{2\pi}-(t-x)^{2}}\right\}  \right]
dx
\]
identical to the southeastern triangle. \ The southern circular cap has area\
\[
2%
%TCIMACRO{\dint \limits_{p_{5}}^{t}}%
%BeginExpansion
{\displaystyle\int\limits_{p_{5}}^{t}}
%EndExpansion
\left[  0-\left\{  t-\sqrt{\dfrac{1}{2\pi}-(t-x)^{2}}\right\}  \right]  dx
\]
identical to the western circular cap. \ The southwestern triangle has area%
\[%
%TCIMACRO{\dint \limits_{0}^{p_{5}}}%
%BeginExpansion
{\displaystyle\int\limits_{0}^{p_{5}}}
%EndExpansion
\left[  \left\{  t-\sqrt{\dfrac{1}{2\pi}-(t-x)^{2}}\right\}  -0\right]  dx
\]
and the northeastern circular cap has area
\[
2%
%TCIMACRO{\dint \limits_{p_{1}}^{1/2}}%
%BeginExpansion
{\displaystyle\int\limits_{p_{1}}^{1/2}}
%EndExpansion
\left[  \left\{  t+\sqrt{\dfrac{1}{2\pi}-(t-x)^{2}}\right\}  -(1-x)\right]
dx+2%
%TCIMACRO{\dint \limits_{1/2}^{p_{7}}}%
%BeginExpansion
{\displaystyle\int\limits_{1/2}^{p_{7}}}
%EndExpansion
\left[  \left\{  t+\sqrt{\dfrac{1}{2\pi}-(t-x)^{2}}\right\}  -x\right]  dx.
\]
Adding these areas and differentiating with respect to $t$, we find that the
best $t$~is%
\[
t=\frac{1}{4}\left(  1+\frac{1}{\pi}\right)  =0.3295774715459476678844418...
\]
corresponding to an asymmetry $\alpha(T_{1})\approx0.450$. \ 

\subsection{$30^{\circ}$-$60^{\circ}$-$90^{\circ}$ Triangle}

Since $\left\vert T_{2}\right\vert =\sqrt{3}/2$, the optimal circle $C$ has
the form%
\[
(x-s)^{2}+(y-t)^{2}=\dfrac{\sqrt{3}}{2\pi}.
\]
The intersection of the line $\sqrt{3}x+y=\sqrt{3}$ and $C$ yields two points
$(p_{1},q_{1})$, $(p_{2},q_{2})$ with $p_{1}<p_{2}$, $q_{1}>q_{2}$. \ The
intersection of the line $x=0$ and $C$ yields two points $(p_{3},q_{3})$,
$(p_{4},q_{4})$ with $p_{3}=p_{4}=0$, $q_{3}>q_{4}$. \ The intersection of the
line $y=0$ and $C$ yields two points $(p_{5},q_{5})$, $(p_{6},q_{6})$ with
$p_{5}<p_{6}$, $q_{5}=q_{6}=0$. \ Everything can be made explicit, for
example,%
\[
p_{1},p_{2}=\frac{1}{4}\left(  3+s-\sqrt{3}t\mp\sqrt{-3+\frac{2\sqrt{3}}{\pi
}+6s-3s^{2}+2\sqrt{3}t-2\sqrt{3}st-t^{2}}\right)  ,
\]
\
\[
q_{1},q_{2}=\frac{1}{4}\left(  \sqrt{3}-\sqrt{3}s+3t\pm\sqrt{-9+\frac
{6\sqrt{3}}{\pi}+18s-9s^{2}+6\sqrt{3}t-6\sqrt{3}st-3t^{2}}\right)  ,
\]%
\[%
\begin{array}
[c]{ccc}%
p_{5}=s-\sqrt{\dfrac{\sqrt{3}}{2\pi}-t^{2}}, &  & q_{4}=t-\sqrt{\dfrac
{\sqrt{3}}{2\pi}-s^{2}}.
\end{array}
\]
The northwestern triangle has area\
\[%
%TCIMACRO{\dint \limits_{0}^{p_{1}}}%
%BeginExpansion
{\displaystyle\int\limits_{0}^{p_{1}}}
%EndExpansion
\left[  \sqrt{3}(1-x)-\left\{  t+\sqrt{\dfrac{\sqrt{3}}{2\pi}-(s-x)^{2}%
}\right\}  \right]  dx
\]
and the southeastern triangle has area%
\[%
%TCIMACRO{\dint \limits_{0}^{q_{2}}}%
%BeginExpansion
{\displaystyle\int\limits_{0}^{q_{2}}}
%EndExpansion
\left[  \frac{1}{\sqrt{3}}\left(  \sqrt{3}-y\right)  -\left\{  s+\sqrt
{\dfrac{\sqrt{3}}{2\pi}-(t-y)^{2}}\right\}  \right]  dy.
\]
\ The southern circular cap has area\
\[
2%
%TCIMACRO{\dint \limits_{p_{5}}^{s}}%
%BeginExpansion
{\displaystyle\int\limits_{p_{5}}^{s}}
%EndExpansion
\left[  0-\left\{  t-\sqrt{\dfrac{\sqrt{3}}{2\pi}-(s-x)^{2}}\right\}  \right]
dx
\]
and the western circular cap has area\
\[
2%
%TCIMACRO{\dint \limits_{q_{4}}^{t}}%
%BeginExpansion
{\displaystyle\int\limits_{q_{4}}^{t}}
%EndExpansion
\left[  0-\left\{  s-\sqrt{\dfrac{\sqrt{3}}{2\pi}-(t-y)^{2}}\right\}  \right]
dy.
\]
The southwestern triangle has area%
\[%
%TCIMACRO{\dint \limits_{0}^{p_{5}}}%
%BeginExpansion
{\displaystyle\int\limits_{0}^{p_{5}}}
%EndExpansion
\left[  \left\{  t-\sqrt{\dfrac{\sqrt{3}}{2\pi}-(s-x)^{2}}\right\}  -0\right]
dx
\]
and the northeastern circular cap has area
\[%
%TCIMACRO{\dint \limits_{q_{2}}^{q_{1}}}%
%BeginExpansion
{\displaystyle\int\limits_{q_{2}}^{q_{1}}}
%EndExpansion
\left[  \left\{  s+\sqrt{\dfrac{\sqrt{3}}{2\pi}-(t-y)^{2}}\right\}  -\frac
{1}{\sqrt{3}}\left(  \sqrt{3}-y\right)  \right]  dy.
\]
Adding these areas and differentiating with respect to both $s$ and $t$, we
find that the best $(s,t)$~is%
\[
s=0.3719164279770188862100673...,
\]%
\[
t=0.4794617554511785131491672...
\]
corresponding to an asymmetry $\alpha(T_{2})\approx0.517$. \ 

\subsection{$6$-$9$-$13$ Triangle}

Since $\left\vert T_{3}\right\vert =4\sqrt{35}$, the optimal circle $C$ has
the form%
\[
(x-s)^{2}+(y-t)^{2}=\dfrac{4\sqrt{35}}{\pi}.
\]
The intersection of the line $4\sqrt{35}x+31y=24\sqrt{35}$ and $C$ yields two
points $(p_{1},q_{1})$, $(p_{2},q_{2})$ with $p_{1}<p_{2}$, $q_{1}>q_{2}$.
\ The intersection of the line $4\sqrt{35}x+13y=0$ and $C$ yields two points;
we select $(p_{3},q_{3})$ to be the point with the larger $y$-coordinate.
\ The intersection of the line $y=0$ and $C$ yields two points $(p_{4},q_{4}%
)$, $(p_{5},q_{5})$ with $p_{4}<p_{5}$, $q_{4}=q_{5}=0$. \ Everything can be
made explicit, for example,%
\[%
\begin{array}
[c]{ccc}%
p_{1},p_{2}=\dfrac{\kappa\mp31\sqrt{\pi\xi}}{1521\pi}, &  & q_{1},q_{2}%
=\dfrac{4\left(  \lambda\pm\sqrt{\frac{35}{\pi}\xi}\right)  }{1521}%
\end{array}
\]
where%
\[
\xi=6084\sqrt{35}-20160\pi+6720\pi s-560\pi s^{2}+1488\sqrt{35}\pi
t-248\sqrt{35}\pi st-961\pi t^{2},
\]%
\[%
\begin{array}
[c]{ccc}%
\kappa=3360\pi+961\pi s-124\sqrt{35}\pi t, &  & \lambda=186\sqrt{35}%
-31\sqrt{35}s+140t;
\end{array}
\]%
\[%
\begin{array}
[c]{ccc}%
p_{3}=\dfrac{13\left(  \mu-\sqrt{\pi\eta}\right)  }{729\pi}, &  & q_{3}%
=\dfrac{4\left(  \nu+\sqrt{\frac{35}{\pi}\eta}\right)  }{729}%
\end{array}
\]
where%
\[
\eta=2916\sqrt{35}-560\pi s^{2}-104\sqrt{35}\pi st-169\pi t^{2},
\]%
\[%
\begin{array}
[c]{ccc}%
\mu=13\pi s-4\sqrt{35}\pi t, &  & \nu=-13\sqrt{35}s+140t;
\end{array}
\]%
\[%
\begin{array}
[c]{c}%
p_{4},p_{5}=s\mp\sqrt{\dfrac{4\sqrt{35}}{\pi}-t^{2}}.
\end{array}
\]

The northwestern triangle has area\
\[%
%TCIMACRO{\dint \limits_{-\frac{13}{3}}^{p_{3}}}%
%BeginExpansion
{\displaystyle\int\limits_{-\frac{13}{3}}^{p_{3}}}
%EndExpansion
\left[  \left\{  -\tfrac{4\sqrt{35}}{31}(x-6)\right\}  -\left\{
-\tfrac{4\sqrt{35}}{13}x\right\}  \right]  dx+%
%TCIMACRO{\dint \limits_{p_{3}}^{p_{1}}}%
%BeginExpansion
{\displaystyle\int\limits_{p_{3}}^{p_{1}}}
%EndExpansion
\left[  \left\{  -\tfrac{4\sqrt{35}}{31}(x-6)\right\}  -\left\{
t+\sqrt{\tfrac{4\sqrt{35}}{\pi}-(s-x)^{2}}\right\}  \right]  dx
\]
and the southeastern triangle has area%
\[%
%TCIMACRO{\dint \limits_{p_{5}}^{p_{2}}}%
%BeginExpansion
{\displaystyle\int\limits_{p_{5}}^{p_{2}}}
%EndExpansion
\left[  \left\{  t-\sqrt{\tfrac{4\sqrt{35}}{\pi}-(s-x)^{2}}\right\}
-0\right]  dx+%
%TCIMACRO{\dint \limits_{p_{2}}^{6}}%
%BeginExpansion
{\displaystyle\int\limits_{p_{2}}^{6}}
%EndExpansion
\left[  \left\{  -\tfrac{4\sqrt{35}}{31}(x-6)\right\}  -0\right]  dx.
\]
The southwestern circular cap has area%
\[
2%
%TCIMACRO{\dint \limits_{p_{4}}^{s}}%
%BeginExpansion
{\displaystyle\int\limits_{p_{4}}^{s}}
%EndExpansion
\left[  0-\left\{  t-\sqrt{\tfrac{4\sqrt{35}}{\pi}-(s-x)^{2}}\right\}
\right]  dx+%
%TCIMACRO{\dint \limits_{0}^{q_{3}}}%
%BeginExpansion
{\displaystyle\int\limits_{0}^{q_{3}}}
%EndExpansion
\left[  \left\{  -\tfrac{13}{4\sqrt{35}}y\right\}  -\left\{  s-\sqrt
{\tfrac{4\sqrt{35}}{\pi}-(t-y)^{2}}\right\}  \right]  dy
\]
and the northeastern circular cap has area%
\[
2%
%TCIMACRO{\dint \limits_{p_{1}}^{s}}%
%BeginExpansion
{\displaystyle\int\limits_{p_{1}}^{s}}
%EndExpansion
\left[  \left\{  t+\sqrt{\tfrac{4\sqrt{35}}{\pi}-(s-x)^{2}}\right\}
-q_{1}\right]  dx+%
%TCIMACRO{\dint \limits_{q_{2}}^{q_{1}}}%
%BeginExpansion
{\displaystyle\int\limits_{q_{2}}^{q_{1}}}
%EndExpansion
\left[  \left\{  s+\sqrt{\tfrac{4\sqrt{35}}{\pi}-(t-y)^{2}}\right\}  -\left\{
6-\tfrac{31}{4\sqrt{35}}y\right\}  \right]  dy.
\]
Adding these areas and differentiating with respect to both $s$ and $t$, we
find that the best $(s,t)$~is%
\[
s=0.9999634051829363409671652...,
\]%
\[
t=2.4097948974186280609774486...
\]
corresponding to an asymmetry $\alpha(T_{3})\approx0.694$.

The fairly arbitrary triangle $T_{3}$, in particular, serves as a benchmark in
\cite{Kb} to distinguish various centers. \ Our preceding value $t$ is the
perpendicular distance from the equiareal disk center to the shortest triangle
side. \ Since the numerical value $2.409...$ does not appear in the database,
we infer that this center is new.

\section{Illuminating Center}

Let $\Lambda$ denote a light source in three-dimensional space of luminosity
$L$. The amount of light an observer $\Theta$ receives from $\Lambda$ is
called its brightness, measured in lumens per unit area. \ Brightness is
inversely proportional to the square of distance between $\Lambda$ and
$\Theta$. \ (Reason: brightness is constant on the sphere of radius $R$,
center $\Lambda$ and thus is equal to $L/\left(  4\pi R^{2}\right)  $).
\ Although we will focus solely on light sources in the \textit{plane}
(streetlights in a triangular park), the preceding \textit{spatial} definition
of brightness is the basis of our model. \ \ 

Given a triangle $T$ (more precisely, its interior), the total brightness%
\[%
%TCIMACRO{\diint \limits_{T}}%
%BeginExpansion
{\displaystyle\iint\limits_{T}}
%EndExpansion
\frac{1}{(x-s)^{2}+(y-t)^{2}}\,dx\,dy
\]
is the quantity we would wish to maximize with respect to $(s,t)\in T$. \ It
turns out that this integral is divergent and
\[%
%TCIMACRO{\diint \limits_{T\smallsetminus D_{\varepsilon}}}%
%BeginExpansion
{\displaystyle\iint\limits_{T\smallsetminus D_{\varepsilon}}}
%EndExpansion
\frac{1}{(x-s)^{2}+(y-t)^{2}}\,dx\,dy
\]
must be examined instead, where $D_{\varepsilon}$ is the disk of radius
$\varepsilon>0$, center $(s,t)$. \ Further, $(s,t)$ cannot be close to
$\partial T$ (the boundary of $T$), that is, we must restrict $(s,t)\in
T\smallsetminus\left(  T\cap N_{\varepsilon}\right)  $ where $N_{\varepsilon}$
is the $\varepsilon$-tubular neighborhood of $\partial T$. \ Under such
conditions, in the limit as $\varepsilon\rightarrow0^{+}$, a geometric
characterization of the maximum point $P=(s,t)$ becomes available. \ 

Let $T$ possess vertices $A$, $B$, $C$. \ Select any two distinct vertices $U
$, $V$ and note that the semiperimeter of subtriangle $UPV$ is
\[
\sigma=\frac{\sqrt{(U-P)\cdot(U-P)}+\sqrt{(V-P)\cdot(V-P)}+\sqrt
{(U-V)\cdot(U-V)}}{2}.
\]
Using Heron's formula, it follows that the ratio of inner angle to area:%
\begin{align*}
\rho(UPV)  &  =\frac{\angle UPV}{\left\vert UPV\right\vert }\\
&  =\tfrac{\arccos\left(  \tfrac{(U-P)\cdot(V-P)}{\sqrt{(U-P)\cdot(U-P)}%
\sqrt{(V-P)\cdot(V-P)}}\right)  }{\sqrt{\left(  \sigma-\sqrt{(U-P)\cdot
(U-P)}\right)  \left(  \sigma-\sqrt{(V-P)\cdot(V-P)}\right)  \left(
\sigma-\sqrt{(U-V)\cdot(U-V)}\right)  \sigma}}%
\end{align*}
satisfies $\rho(APB)=\rho(BPC)=\rho(CPA)$, a remarkable fact \cite{Sh}! \ For
$T_{1}$, we solve the resulting equations, obtaining
\[
s=t=0.3082756986146550422567206...;
\]
for $T_{2}$,%
\[
s=0.3516876887676632055410277...,
\]%
\[
t=0.4491286165669552235961426...;
\]
and, for $T_{3}$,%
\[
s=0.8345011650594754190821304...,
\]%
\[
t=2.0031487728161056257679347....
\]
\ Our preceding value $t$ is the perpendicular distance from the illuminating
center to the shortest triangle side. \ Since the numerical value $2.003...$
does not appear in the ETC\ database, we infer that this center is new.

Let us revisit the definition of brightness. \ Had a \textit{planar}
definition been adopted -- measured in lumens per unit length -- then
brightness would be inversely proportional to the distance itself between
$\Lambda$ and $\Theta$. \ (Reason: brightness would be constant on the circle
of radius $R$, center $\Lambda$ and thus would be equal to $L/\left(  2\pi
R\right)  $). \ This scenario yields exactly the same formulation as that
underlying the electrostatic center \cite{AK}. \ \ 

\section{Thermodynamic Center}

Here, the triangle $T$ (more precisely, its interior) is assumed to be a heat
conductor with initial temperature $1$ while its boundary $\partial T$ is held
at temperature $0$ always. \ Heat $u$ will dissipate as time $t$ increases
according to the following initial-boundary value problem:
\[%
\begin{array}
[c]{ccc}%
\dfrac{\partial u}{\partial t}=\dfrac{\partial^{2}u}{\partial x^{2}}%
+\dfrac{\partial^{2}u}{\partial y^{2}} &  & \text{for }(x,y,t)\in
T\times(0,\infty),
\end{array}
\]%
\[%
\begin{array}
[c]{ccc}%
u(x,y,t)=1 &  & \text{for }(x,y,t)\in T\times\{0\},
\end{array}
\]%
\[%
\begin{array}
[c]{ccc}%
u(x,y,t)=0 &  & \text{for }(x,y,t)\in\partial T\times(0,\infty)
\end{array}
\]
but does so non-uniformly: its density gathers around a unique maximum point
$(x_{\infty},y_{\infty})$ as $t\rightarrow\infty$ \cite{Pa, BMS}. The point
$(x_{\infty},y_{\infty})$ is, in fact, the unique extreme point of the first
Laplacian eigenfunction for $T$. \ The eigenfunction for $T_{1}$ is \cite{Mc,
Si, DP, F1, F2, F3}%
\[
\sin(\pi x)\sin(2\pi y)+\sin(2\pi x)\sin(\pi y)
\]
with%
\[
x_{\infty}=y_{\infty}=\frac{1}{\pi}\operatorname{arcsec}\left(  \sqrt
{3}\right)  =0.3040867239846963649145722...
\]
and the eigenfunction for $T_{2}$ is%
\[
\sin\left(  \frac{\pi x}{3}\right)  \sin\left(  \sqrt{3}\pi y\right)
+\sin\left(  \frac{4\pi x}{3}\right)  \sin\left(  \frac{2\pi y}{\sqrt{3}%
}\right)  +\sin\left(  \frac{5\pi x}{3}\right)  \sin\left(  \frac{\pi y}%
{\sqrt{3}}\right)
\]
with%
\[
x_{\infty}=0.3558473606263811208579681...,
\]%
\[
y_{\infty}=0.4255359610370576630888604....
\]
It remains to compute the eigenfunction for $T_{3}$, but no closed-form
expression for this exists. \ A\ numerical computation using the Matlab
\textit{pdeeig} tool yields $(x_{\infty},y_{\infty})\approx(0.88,1.91)$. \ See
Figure 4. \ Significantly higher precision will be needed to ascertain whether
$y_{\infty}$ appears in the ETC\ database; we are hopeful that techniques in
\cite{HT} might save the day.%
%TCIMACRO{\FRAME{ftbpFU}{5.828in}{2.7518in}{0pt}{\Qcb{First Laplacian
%eigenfunction for $T_{3}$ using Matlab's \QTR{it}{pdesurf}}}{}{tri4.eps}%
%width 5.828in;  height 2.7518in;  depth 0pt;  original-width 17.7849in;
%original-height 8.3679in;  cropleft "0";  croptop "1";  cropright "1";
%cropbottom "0";  filename 'tri4.eps';file-properties "XNPEU";}} }%
%BeginExpansion
\begin{figure}[ptb]%
\centering
\includegraphics[
height=2.7518in,
width=5.828in
]%
{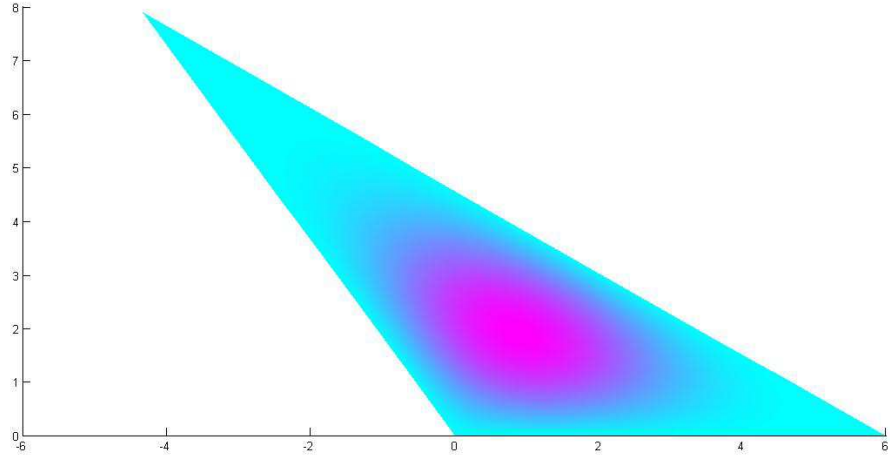}%
\caption{First Laplacian eigenfunction for $T_{3}$ using Matlab's
\textit{pdesurf}}%
\end{figure}
%EndExpansion

\section{Acknowledgements}

I am thankful to Harri Hakula for refining my preceding estimate to
$(x_{\infty},y_{\infty})=(0.88047...,1.91599...)$ via Mathematica \cite{HH},
Stefan Steinerberger for introducing me to Fraenkel asymmetry, and the authors
of \cite{AK, BM} for writing inspirational papers.

\end{document}